\documentclass[11pt]{article}
\newcommand{\QSL}{\mathit{QSL}}
\newcommand{\PF}{\mathrm{PF}}
\newcommand{\UPF}{\mathrm{UPF}}
\newcommand{\Rep}{\mathrm{Rep}}
\newcommand{\Min}{\operatorname{Min}}
\newcommand{\MIN}{\operatorname{MIN}}
\newcommand{\COL}{\operatorname{COL}}
\newcommand{\Max}{\operatorname{Max}}
\newcommand{\MAX}{\operatorname{MAX}}
\newcommand{\ROW}{\operatorname{ROW}}
\newcommand{\MHS}{\mathcal{M}_\mathrm{HS}}
\newcommand{\Mcb}{\mathcal{M}_\mathrm{cb}}
\newcommand{\KG}{K_{\mathbb G}}
\renewcommand{\baselinestretch}{1.2}
\newcommand{\dated}{\mbox{} \hfill {\small [{\tt \today}]}} 
\usepackage{amsthm,enumerate}
\theoremstyle{plain}
\newtheorem{theorem}{Theorem}[section]
\newtheorem{lemma}[theorem]{Lemma}
\newtheorem{corollary}[theorem]{Corollary}
\newtheorem{proposition}[theorem]{Proposition}
\theoremstyle{definition}
\newtheorem{definition}[theorem]{Definition}
\theoremstyle{remark}
\newtheorem*{remark}{Remark}

\newtheorem*{rems}{Remarks}
\newtheorem*{exs}{Examples}
\newenvironment{remarks}{\begin{rems}\begin{enumerate}}{\end{enumerate}\end{rems}}

\newenvironment{items}{\begin{enumerate}[\rm (i)]}{\end{enumerate}}
\newenvironment{alphitems}{\begin{enumerate}[\rm (a)]}{\end{enumerate}}

 \usepackage{amsmath,amssymb,amsfonts,diagrams}
\newenvironment{keywords}{\noindent\small {\it Keywords\/}:}{\vskip 4pt}
\newenvironment{classification}{\noindent\small 2000 {\it Mathematics Subject
Classification\/}:}{\vskip 12pt}

\newcommand{\comps}{{\mathbb C}}

\newcommand{\posints}{{\mathbb N}}

\newcommand{\tensor}{\otimes}
\newcommand{\ttensor}{\tilde{\otimes}}
\newcommand{\Tensor}{\hat{\otimes}}

\newcommand{\cstar}{{C^\ast}}

\newcommand{\id}{{\mathrm{id}}}
\newcommand{\cb}{{\mathrm{cb}}}

\newcommand{\tr}{{\operatorname{tr}}}

\newcommand{\A}{{\mathfrak A}}

\newcommand{\Hilbert}{{\mathfrak H}}

\newcommand{\CB}{{\cal CB}}
\newcommand{\op}{{\mathrm{op}}}

\newcommand{\VN}{\operatorname{VN}}
\newcommand{\PM}{\operatorname{PM}}

\title{Column and row operator spaces over $\QSL_p$-spaces \\ and their use in abstract harmonic analysis}
\author{\textit{Matthias Neufang} \and \textit{Volker Runde}}
\date{}
\begin{document}
\maketitle
\begin{abstract}
The notions of column and row operator space were extended by A.\ Lambert from Hilbert spaces to general Banach spaces. In this paper, we use column and row spaces over quotients of subspaces of general $L_p$-spaces to equip several Banach algebras occurring naturally in abstract harmonic analysis with canonical, yet not obvious operator space structures that turn them into completely bounded Banach algebras. We use these operator space structures to gain new insights on those algebras.
\end{abstract}
\section*{Introduction}
The Fourier algebra $A(G)$ of a general locally compact group $G$ was introduced by P.\ Eymard in \cite{Eym}. If $G$ is abelian with dual group $\hat{G}$, then $A(G)$ is just $L_1(\hat{G})$ via the Fourier transform. As the predual of the group von Neumann algebra, $A(G)$ has a canonical structure as an abstract operator space (see \cite{ER}, \cite{Pau}, or \cite{Pis} for the theory of operator spaces), turning it into a completely contractive Banach algebra. In the past decade and a half, operator space theoretic methods have given new momentum to the study of $A(G)$ (see \cite{IS}, \cite{NRS}, or \cite{Rua}, for example), yielding new insights, even if the problem in question seemed to have nothing to do with operator spaces (\cite{FKLS} or \cite{FR}).
\par 
The definition of $A(G)$ can be extended to an $L_p$-context: instead of restricting oneself the left regular representation of $G$ on $L_2(G)$, one considers the left regular representation of $G$ on $L_p(G)$ for general $p \in (1,\infty)$. This approach leads to the \emph{Fig\`a-Talamanca--Herz algebras} $A_p(G)$, which were introduced by C.\ Herz in \cite{Her1} and further studied in \cite{Her2}. Ever since, the Fig\`a-Talamanca--Herz algebras have been objects of independent interest in abstract harmonic analysis. At the first glance, it may seem that the passage from $L_2(G)$ to $L_p(G)$ for $p \neq 2$ is of little significance, and, indeed, many (mostly elementary) properties of $A(G)$ can be established for $A_p(G)$ with $p \neq 2$ along the same lines. However, the lack of von Neumann algebraic methods for operator algebras on $L_p$-spaces for $p \neq 2$ has left other problems, which have long been solved for $A(G)$, wide open for $A_p(G)$. For instance, any closed subgroup of $G$ is a set of synthesis for $A(G)$ (\cite{TT}) whereas the corresponding statement for $A_p(G)$ with $p \neq 2$ is still wide open.
\par 
As the Fig\`a-Talamanca--Herz algebras have no obvious connections with operator algebras on Hilbert space, it appears at first glance that operator space theoretic methods are of very limited use when dealing with $A_p(G)$ for $p \neq 2$. There is a notion of $p$-completely boundedness for general $p \in (1,\infty)$ with $2$-complete boundedness just being usual complete boundedness, and an abstract theory based on $p$-complete boundedness can be developed---called $p$-operator space theory in \cite{Daw}---that parallels operator space theory (\cite{LeM}). There are indeed applications of $p$-complete boundedness to Fig\`a-Talamanca--Herz algebras (see \cite{Fen} and \cite{Daw}). Alas, as pointed out in \cite{Daw}, there is no suitable Hahn--Banach theorem for $p$-completely bounded maps, so that the duality theory of $p$-operator spaces has to be fairly limited.
\par 
In \cite{LNR}, A.\ Lambert and the authors pursued a different approach to putting operator spaces to work on Fig\`a-Talamanca--Herz algebras. In his doctoral thesis \cite{Lam}, Lambert extended the notions of column and row operator space, which are canonical over Hilbert space, to general Banach spaces. This allows, for $p \in (1,\infty)$, to equip ${\cal B}(L_p(G))$ for any $p \in (1,\infty)$ with an operator space structure, which, for $p = 2$, is the canonical one. This, in turn, can be used to equip $A_p(G)$---for any $p \in (1,\infty)$---with an operator space structure in the usual sense, making it a completely bounded Banach algebra. With respect to this operator space structure, \cite[Theorem 3.6]{Rua} extends to Fig\`a-Talamanca--Herz algebras: $G$ is amenable if and only if $A_p(G)$ is operator amenable for one---and, equivalently, all---$p \in (1,\infty)$.
\par 
In the present paper, we continue the work begun in \cite{LNR} and link it with the paper \cite{RunBp} by the second author. Most of it is devoted to extending operator space theoretic results known to hold for the Fourier algebra and (reduced) Fourier--Stieltjes algebra of a locally compact group to the suitable generalizations in a general $L_p$-context. In particular, we show that, for any $p \in (1,\infty)$, the Banach algebra $B_p(G)$ introduced in \cite{RunBp} can be turned into a completely bounded Banach algebra in a canonical manner, and we obtain an $L_p$-generalization of \cite[Theorem 4.4]{RS}.
\section{Preliminaries}
In this section, we recall some of background from \cite{LNR} and \cite{RunBp}. We shall throughout rely heavily on those papers, and the reader is advised to have them at hand.
\subsection{Column and row operators spaces over Banach spaces}
The notions of column and row operator space of Hilbert space are standard in operator space theory (\cite[3.4]{ER}). In \cite{Lam}, Lambert extended these notions to general Banach spaces. As his construction is fairly involved, we will only sketch it very briefly here and refer to \cite[Sections 2 and 3]{LNR} instead (and to \cite{Lam} for more details). Throughout the paper, we adopt the notation from \cite{LNR}.
\par 
Lambert introduces a category---called \emph{operator sequence spaces}---that can be viewed as an intermediary between Banach spaces and operator spaces, and defines functors
\[
  \min, \max \!: \{\text{Banach spaces} \} \to \{ \text{operator sequence spaces} \}
\]
and
\[
  \Min, \Max \!: \{\text{operator sequence spaces} \} \to \{ \text{operator spaces} \}
\]
such that $\Min \circ \min = \MIN$ and $\Max \circ \max = \MAX$. He then defines
\[
  \COL, \ROW \!: \{\text{Banach spaces} \} \to \{ \text{operator spaces} \}
\]
as
\[
  \COL := \Min \circ \max \qquad\text{and}\qquad \ROW := \Max \circ \min.
\]
For any Banach space $E$, the operator spaces $\COL(E)$ and $\ROW(E)$ are homogeneous and satisfy 
\[
  \COL(E)^\ast = \ROW(E^\ast) \qquad\text{and}\qquad \ROW(E)^\ast = \COL(E^\ast).
\]
By \cite{Math}, these definitions coincide with the usual ones in the case of a Hilbert space.
\par 
We would like to point out that column and row operator spaces in a general Banach space context display a behavior quite different from in the Hilbert space setting, as shown by the following two propositions:
\begin{proposition} \label{colmax}
For any subhomogeneous $\cstar$-algebra $\A$, we have canonical completely bounded isomorphisms
\[
  \COL(\A) \cong \MIN(\A) \qquad\text{and}\qquad \ROW(\A) \cong \MAX(\A).
\]
\end{proposition}
\begin{proof}
Let $\A$ be a subhomogeneous $\cstar$-algebra. As $\id_\A \!: \min(\A) \to \max(\A)$ is sequentially bounded by \cite[Satz 2.2.25]{Lam}, 
\[
  \id_\A \!: \MIN(\A) = \Min(\min(\A)) \to \Min(\max(\A)) = \COL(\A)
\]
and
\[
  \id_\A \!: \ROW(\A) = \Max(\min(\A)) \to \Max(\max(\A)) = \MAX(\A)
\]
are completely bounded.
\end{proof}
\begin{remark}
Let $\Hilbert$ be an infinite-dimensional Hilbert space, and let $\iota \!: \Hilbert \to {\cal C}(\Omega)$ be an isometric embedding into the continuous functions onto some compact Hausdorff space $\Omega$. As $\COL({\cal C}(\Omega))\cong \MIN({\cal C}(\Omega))$ by Proposition \ref{colmax} whereas $\id_\Hilbert \!: \MIN(\Hilbert) \to \COL(\Hilbert)$ is not completely bounded, this means that $\COL$ does not respect subspaces; in a similar way, we see that $\COL$ does not respect quotients either. By duality, the same is true for $\ROW$.
\end{remark}
\par
Given an operator space $E$, we denote its \emph{opposite operator space} by $E^\op$ (see \cite[pp.\ 43--44]{Pis0}). It is immediate that $\COL(\Hilbert)^\op = \ROW(\Hilbert)$ and $\ROW(\Hilbert)^\op = \COL(\Hilbert)$ for any Hilbert space $\Hilbert$. For a general Banach space $E$, we have:
\begin{proposition} \label{opposite}
For any Banach space $E$, the identity on $E$ induces a complete contraction from $\COL(E)$ to $\ROW(E)^\op$, which, in general, fails to have a completely bounded inverse.
\end{proposition}
\begin{proof}
From the definition of $\ROW$, it is obvious that $\id_E \!: \min(E) \to C(\ROW(E)^\op)$ is a sequential isometry, so that
\[
  \id_E \!: \COL(E) = \Max(\min(E)) \to \ROW(E)^\op
\]
is a complete contraction by \cite[Satz 4.1.12]{Lam}.
\par 
On the other hand, if $E$ is a subhomogeneous $\cstar$-algebra, then we have $\ROW(\A)^\op \cong \MAX(\A)^\op = \MAX(\A)$ by Proposition \ref{colmax} whereas $\COL(\A) \cong \MIN(\A)$. Hence, $\id \!: \ROW(\A)^\op \to \COL(\A)$ cannot be completely bounded unless $\dim \A < \infty$ (\cite[Theorem 14.3(iii)]{Pau}).
\end{proof}
\subsection{Representations on $\QSL_p$-spaces}
By a \emph{representation} of a locally compact group $G$ on a Banach space, we mean a pair $(\pi,E)$, where $E$ is a Banach space and $\pi$ is homomorphism from $G$ into the group of invertible isometries on $E$ which is continuous with respect to the given topology on $G$ and the strong operator topology on ${\cal B}(E)$. (This is somewhat more restrictive than the usual use of the term; as we will not consider any other kind of representation, however, we prefer to keep our terminology short.) 
\par 
Given two representations $(\pi,E)$ and $(\rho,F)$ of a locally compact group $G$, we 
\begin{itemize}
\item call $(\pi,E)$ and $(\rho,F)$ \emph{equivalent} if there is an invertible isometry $V \!: E \to F$ with
\[
  V\pi(x) V^{-1} = \rho(x) \qquad (x \in G),
\]
\item call $(\rho,F)$ a \emph{subrepresentation} of $(\pi,E)$ if $F$ is a closed subspace of $E$ and $\rho(x) = \pi(x) |_F$ holds for all $x \in F$, and
\item say that $(\rho,F)$ is \emph{contained} in $(\pi,E)$ if $(\rho,F)$ is equivalent to a subrepresentation of $(\pi,E)$, in which case, we write $(\rho,F) \subset (\pi,E)$.
\end{itemize}
\par 
Throughout, we shall often identify a particular representation with its equivalence class in order to avoid pedantry.
\par 
Given a locally compact group $G$ and a representation $(\pi,E)$ of $G$, we obtain a representation of the group algebra $L_1(G)$ on $E$, i.e., a contractive algebra homomorphism from $L_1(G)$ into ${\cal B}(E)$, which we denote by $\pi$ as well, via
\begin{equation} \label{integ}
  \pi(f) := \int_G f(x) \pi(x) \, dx \qquad (f \in L_1(G)),
\end{equation}
where the integral converges in the strong operator topology. Conversely, if $\pi \!: L_1(G) \to {\cal B}(E)$ is a representation that is non-degenerate, i.e., the span of $\{ \pi(f)\xi : f \in L_1(G), \, \xi \in E \}$ is dense in $E$, then it arises from a representation of $G$ on $E$ via (\ref{integ}).
\par 
In this paper, we are interested in representations on $\QSL_p$-spaces with $p \in (1,\infty)$, i.e., on Banach spaces that are isometrically isomorphic to quotients of subspaces---or, equivalently, subspace of quotients---of the usual $L_p$-spaces. By \cite[{\S}4, Theorem 2]{Kwa}, these are precisely the $p$-spaces of \cite{Her2}. For a locally compact group $G$ and $p \in (1,\infty)$, we denote by $\Rep_p(G)$ the collection of all (equivalence classes of) representations of $G$ on a $\QSL_p$-space. 
\par 
For the following definition, recall that, for $p \in (1,\infty)$, any $\QSL_p$-space $E$ is reflexive, so that ${\cal B}(E)$ is a dual Banach space in a canonical way, so that we can speak of a weak$^\ast$ topology.
\begin{definition}
Let $G$ be a locally compact group, let $p \in (1,\infty)$, and let $(\pi,E) \in \Rep_p(G)$. Then:
\begin{alphitems}
\item the algebra $\PF_{p,\pi}(G)$ of \emph{$p$-pseudofunctions associated with $(\pi,E)$} is the norm closure of $\pi(L_1(G))$ in ${\cal B}(E)$;
\item the algebra $\PM_{p,\pi}(G)$ of \emph{$p$-pseudomeasures associated with $(\pi,E)$} is the weak$^\ast$ closure of $\pi(L_1(G))$ in ${\cal B}(E)$.
\end{alphitems}
\end{definition}
\par 
If $(\pi,E) = (\lambda_p,L_p(G))$, i.e., the left regular representation of $G$ on $L_p(G)$, we simply speak of $p$-pseudofunctions and $p$-pseudomeasures, as is standard usage, and write $\PF_p(G)$ and $\PM_p(G)$, respectively. 
\section{$\PF_p(G)$ is not an operator algebra for $p \neq 2$}
Let $G$ be a locally compact group. Then $\PF_2(G)$ and $\PM_2(G)$ are the reduced group $\cstar$-algebra $C^\ast_r(G)$ and the group von Neumann algebra $\VN(G)$, respectively. As $\cstar$-subalgebras of ${\cal B}(L^2(G))$, they are operator spaces in a canonical manner. For any Hilbert space $\Hilbert$, the operator spaces ${\cal B}(\Hilbert)$ and $\CB(\COL(\Hilbert))$ are completely isometrically isomorphic. 
\par 
We thus define:
\begin{definition} \label{osdef}
Let $G$ be a locally compact group, let $p \in (1,\infty)$, and let $(\pi,E) \in \Rep_p(G)$. Then the canonical operator space structure of $\PF_{p,\pi}(G)$ and $\PM_{p,\pi}(G)$, respectively, is the one inherited as a subspace of $\CB(\COL(E))$.
\end{definition}
\par 
As $\CB(E)$ is a completely contractive Banach algebra, for every operator space $E$, Definition \ref{osdef} turns $\PF_{p,\pi}(G)$ and $\PM_{p,\pi}(G)$ into completely contractive Banach algebras.
\par 
By an \emph{operator algebra}, we mean a completely contractive Banach algebra that is completely isometrically isomorphic to a---not necessarily self-adjoint---closed subalgebra of ${\cal B}(\Hilbert)$ for some Hilbert space $\Hilbert$. Not every completely bounded Banach algebra is an operator algebra: for instance, $\CB(E)$ for an operator space $E$, is an operator algebra if and only if $E = \COL(\Hilbert)$ for some Hilbert space $E$ (\cite[Theorem 3.4]{Ble}). Nevertheless, there are operator algebras that arise naturally from $L_p$-spaces for $p \neq 1$ (\cite{BLeM}). Hence, the following theorem is somewhat less self-evident than one might think at the first glance:
\begin{theorem} \label{noopalg}
Let $p \in (1,\infty) \setminus \{ 2 \}$, and let $G$ be an amenable, locally compact group containing an infinite abelian subgroup. Then $\PF_p(G)$ and $\PM_p(G)$ are not operator algebras.
\end{theorem}
\par
Before we delve into the proof, we establish some more notation and definitions. 
\par 
For a locally compact group $G$ and $p \in (1,\infty)$, we denote its \emph{Fig\`a-Talamanca--Herz algebra} by $A_p(G)$ (see \cite{Her1}, \cite{Her2}, and \cite{Eym2}). We have a canonical duality $A_p(G)^\ast = \PM_{p'}(G)$, where $p' \in (1,\infty)$ is such that $\frac{1}{p} + \frac{1}{p'} = 1$. Of course, $A_2(G)$ is just Eymard's Fourier algebra $A(G)$. By a \emph{multiplier} of $A_p(G)$, we mean a continuous function $f$ such that $fg \in A_p(G)$ for all $g \in A_p(G)$; for any multiplier $f$ of $A_p(G)$, multiplication induces a linear map $M_f \!: A_p(G) \to A_p(G)$, which is easily seen to be bounded by the closed graph theorem. As $A_p(G)$ is a closed subspace of $\PM_{p'}(G)^\ast$, it inherits a canonical operator space structure via Definition \ref{osdef} (this is the same operator space structure as considered in \cite{LNR}). We thus define the \emph{completely bounded multipliers} of $A_p(G)$ as
\[
  \Mcb(A_p(G)) := \{ f : \text{$f$ is a multiplier of $A_p(G)$ such that $M_f \in \CB(A_p(G))$} \}.
\]
It is easy to see that $\Mcb(A_p(G))$ is a closed subalgebra of $\CB(A_p(G))$. In \cite{LNR}, it was shown that multiplication in $A_p(G)$ is completely bounded, even though not necessarily completely contractive, so that we have a canonical completely bounded inclusion $A_p(G) \subset \Mcb(A_p(G))$.
\par 
We start the proof of Theorem \ref{noopalg} with a lemma:
\begin{lemma} \label{nol1}
Let $G$ be a locally compact group, and let $p \in (1,\infty)$ be such that $\PF_p(G)$ is an operator algebra. Then $A_p(G) \subset \Mcb(A(G))$ holds.
\end{lemma}
\begin{proof}
Let $f \in A_{p'}(G)$, so that $M_f \in \Mcb(A_{p'}(G))$ and, consequently, $M_f^\ast \in \CB(\PM_p(G))$. It is easy to see that $M_f^\ast$ leaves $\PF_p(G)$ invariant, so that we can view $M_f^\ast$ as an element of $\CB(\PF_p(G))$.
\par
Assume that there is a Hilbert space $\Hilbert$ and a completely isometric algebra homomorphism $\iota \!: \PF_p(G) \to {\cal B}(\Hilbert)$. Then $\lambda_p^\# := \iota \circ \lambda_p \!: L_1(G) \to {\cal B}(\Hilbert)$ is a contractive representation of $L_1(G)$ on $\Hilbert$. Let $( e_\alpha )_\alpha$ be an approximate identity for $L_1(G)$ bounded by one, and let $P \in {\cal B}(\Hilbert)$ a weak$^\ast$accumulation point of $(\lambda^\#_p(e_\alpha) )_\alpha$. Then $P$ is a norm one idempotent and thus
a projection. Replacing $\Hilbert$ by $P\Hilbert$, we can therefore suppose that $\lambda^\#_p \!: L_1(G) \to {\cal B}(\Hilbert)$ is a non-degenerate representation of $L_1(G)$. Consequently, it arises from a unitary representation of $G$ on $\Hilbert$, which we denote likewise by $\lambda_p^\#$. 
\par
We can view $M_f^\ast$ as a completely bounded map from $\PF_p(G)$ to ${\cal B}(\Hilbert)$. By the Arveson--Wittstock--Hahn--Banach theorem  (\cite[Theorem 4.1.5]{ER}), we can extend $M_f^\ast$ to a completely bounded map $\widetilde{M^\ast_f} 
\!: {\cal B}(\Hilbert) \to {\cal B}(\Hilbert)$. By \cite[Theorem 8.4]{Pau}, there is another Hilbert space $\mathfrak K$ along with a unital $^\ast$-homomorphism $\pi \!: {\cal B}(\Hilbert) \to {\cal B}(\mathfrak{K})$ and bounded operators $V,W \!: \Hilbert \to \mathfrak K$ such that
\[
  \widetilde{M^\ast_f}(T) = V^\ast \pi(T) W \qquad (T \in {\cal B}(\Hilbert))
\]
and thus
\begin{equation} \label{factor1}
  M^\ast_f(\lambda^\#_p(g)) = V^\ast \pi(\lambda^\#_p(g)) W
  \qquad (g \in L_1(G)).
\end{equation}
Set $\sigma := \pi \circ \lambda^\#_p$, so that $\sigma$ is a contractive representation of $L_1(G)$ on the Hilbert space $\mathfrak K$.
As before for $\lambda_p^\#$, we see that $\sigma$ is non-degenerate and thus is induced via integration by a unitary representation of $G$ on $\mathfrak K$, which we denote likewise by $\sigma$. If $( e_\alpha )_\alpha$ denotes again a bounded approximate identity for $L_1(G)$, we obtain 
\begin{equation} \label{factor2}
  \begin{split}
  f(x) \lambda_p^\#(x) & = M_f^\ast(\lambda_p^\#(x)) \\
  & = \lim_\alpha M_f^\ast(\lambda_p^\#(\delta_x \ast e_\alpha)) \\
  & = \lim_\alpha V^\ast \sigma(\delta_x \ast e_\alpha) W, \qquad\text{by (\ref{factor1})}, \\
  & = V^\ast \sigma(x) W \qquad (x \in G),
  \end{split}
\end{equation}
where the limits are taken in the strong operator topology.
\par
Fix a unit vector $\xi \in \Hilbert$, and define $L, R \!: G \to \Hilbert$ by letting
\[
  L(x) := \sigma(x^{-1})V\lambda_p^\#(x)\xi 
  \quad\text{and}\quad
  R(x) := \sigma(x^{-1}) W \lambda_p^\#(x)\xi
  \qquad (x \in G).
\]
It follows that
\[
  \begin{split}
  \langle L(x),R(y) \rangle & = 
  \langle \sigma(x^{-1})V\lambda_p^\#(x)\xi,
          \sigma(y^{-1})W\lambda_p^\#(y)\xi \rangle \\
  & = \langle \lambda_p^\#(x)\xi, 
          V^\ast \sigma(xy^{-1})W\lambda_p^\#(y)\xi \rangle \\
  & = \langle f(xy^{-1}) \lambda_p^\#(x)\xi, 
          \lambda_p^\#(xy^{-1})\lambda_p^\#(y)\xi \rangle, \qquad
          \text{by (\ref{factor2})}, \\
  & = f(xy^{-1}) \qquad (x,y \in G).
  \end{split}
\]
\par
By the characterization of completely bounded multipliers of $A(G)$ from \cite{BF}---see \cite{Jol} for an alternative proof---, this means that $f \in \Mcb(A(G))$.
\end{proof}
\begin{remark}
Our proof of Lemma \ref{nol1} is patterned after that of the main result of \cite{Jol}.
\end{remark}
\par 
Our next proposition seems to be folklore (and likely true for all infinite $G$), but for lack of a suitable reference, we give a proof that was indicated to us by A.\ Derighetti:
\begin{proposition} \label{nop}
Let $G$ be a locally compact group that contains an infinite abelian subgroup, and let $p \in (1,\infty) \setminus \{ 2 \}$. Then $A_p(G)$ is not contained in $A(G)$.
\end{proposition}
\begin{proof}
Let $H$ be an infinite abelian subgroup of $G$, which we can suppose to be closed, and note that $\PM_p(H) \subsetneq \PM_2(H)$ by \cite[Theorems 4.5.4 and 4.5.5]{Lar}, so that, by duality, $A(H) \subsetneq A_p(H)$.
\par 
Assume that $A_p(G) \subset A(G)$, and let $Q_{p,H} \!: A_p(G) \to A_p(H)$ and $Q_H \!: A(G) \to A(H)$ denote the restriction maps. Since both $Q_{p,H}$ and $Q_H$ are surjective (\cite[Theorem 1b]{Her2}), $A_p(G) \subset A(G)$ thus yields $A_p(H) \subset A(H)$, which contradicts $A(H) \subsetneq A_p(H)$.
\end{proof}
\begin{remark}
Groups to which Proposition \ref{nop} applies are, in particular, all infinite, compact groups (\cite[Theorem 2]{Zel}), all non-metrizable groups (by \cite[(8.7) Theorem]{HR} and \cite[Theorem 2]{Zel} combined), and all connected groups (\cite[Theorem in 4.13]{MZ}).
\end{remark}
\begin{proof}[Proof of Theorem \emph{\ref{noopalg}}]
Since $\PF_p(G)$ is a closed subalgebra of $\PM_p(G)$, it is enough to prove the claim for $\PF_p(G)$.
\par 
By Lemma \ref{nol1}, we have $A_{p'}(G) \subset \Mcb(A(G))$. Since $G$ is amenable, $A(G)$ has a bounded approximate identity by
Leptin's theorem (\cite[Theorem 10.4]{Pie}). Moreover, $A_{p'(G)}$ is a Banach $A(G)$-module by \cite{Her1}, and the bounded approximate identity for $A(G)$ constructed in the proof of \cite[Theorem 10.4]{Pie} is easily seen to also be a bounded approximate identity for the Banach $A(G)$-module. From one of the versions of Cohen's factorization theorem---\cite[Corollary 2.9.26]{Dal}, for instance---, it then follows that each function in $A_{p'}(G)$ is a product of another function in $A_{p'}(G)$ with a function
from $A(G)$. By Lemma \ref{nol1}, $A_{p'}(G) \subset \Mcb(A(G))$ holds, so that $A_{p'}(G) \subset A(G)$. This contradicts Proposition \ref{nop}.
\end{proof}
\section{$\PF_{p'}(G)$ and $B_p(G)$ as completely bounded Banach algebras}
Let $G$ be a locally compact group. Then $A(G) = \VN(G)_\ast$, $B_r(G) = C^\ast_r(G)^\ast$---the \emph{reduced Fourier--Stieltjes algebra} of $G$, and the \emph{Fourier--Stieltjes algebra} $B(G) = \cstar(G)^\ast$, where $\cstar(G)$ denotes the full group $\cstar$-algebra, all have canonical operator space structures turning them into completely contractive Banach algebras. 
\par 
For $p \in (1,\infty)$, the embedding $A_p(G) \subset \PM_{p'}(G)^\ast$ turns $A_p(G)$ into a \emph{completely bounded Banach algebra}, i.e., turns it into an operator space such that multiplication is completely bounded, albeit not necessarily completely contractive (see \cite{LNR} for details).
\par 
For any $p \in (1,\infty)$, the space $\PF_{p'}(G)^\ast$ consists of continuous functions on $G$ and is a Banach algebra under pointwise multiplication (see \cite{Her3} and \cite{Cow}). Moreover, in \cite{RunBp}, the second author defined a unital, commutative Banach algebra $B_p(G)$ containing $\PF_{p'}(G)^\ast$ (\cite[Theorem 6.6(i)]{RunBp}), which, for $p = 2$, is just $B(G)$. In this section, we will adapt the construction from \cite{LNR} to equip both $\PF_{p'}(G)$ and $B_p(G)$ with canonical operator space structures---generalizing those for $B_r(G)$ and $B(G)$ in the $p=2$ case---such that they become completely bounded Banach algebras. 
\par 
We begin the following proposition:
\begin{proposition} \label{tprop}
Let $p \in (1, \infty)$, and let $E$ and $F$ be $\QSL_p$-spaces. Then there is a norm $\| \cdot \|_p$ on the algebraic tensor product $E \tensor F$ with the following properties:
\begin{items}
\item $\| \cdot \|_p$ is a cross norm dominating the injective tensor norm;
\item the completion $E \ttensor_p F$ of $(E \tensor F, \| \cdot \|_p)$ is a $\QSL_p$-space;
\item if $G$ is a locally compact group with $(\pi,E), (\rho, F) \in \Rep_p(G)$, then $(\pi \tensor \rho, E \ttensor_p F) \in \Rep_p(G)$;
\item the bilinear maps
\[
  \COL(E) \times \COL(F) \to \COL(E \ttensor_p F), \quad (\xi,\eta) \mapsto \xi \tensor \eta
\]
and
\[
  \ROW(E) \times \ROW(F) \to \ROW(E \ttensor_p F), \quad (\xi,\eta) \mapsto \xi \tensor \eta
\]
are completely bounded with $\cb$-norm at most $\KG$, the \emph{complex Grothendieck constant}.
\end{items}
Moreover, if $E = L^p(X)$ for some measure space $X$, we can choose $\| \cdot \|_p$ as the norm $L^p(X) \tensor F$ inherits as a subspace of $L^p(X,F)$.
\end{proposition}
\begin{proof}
(i), (ii), and (iii) just summarize \cite[Theorem 3.1]{RunBp} and the ``moreover'' part is clear from an inspection of the proof of that theorem.
\par 
(iv) follows from \cite[Theorem 5.5 and 5.8]{LNR} and the construction of $\| \cdot \|_p$ in \cite{RunBp}.
\end{proof}
\par 
Given a locally compact group $G$, $p \in (1,\infty)$, and $(\pi,E) \in \Rep_{p'}(G)$, let $\PM_{p,,\pi}(G)_\ast$ denote the canonical predual of $\PM_{p',\pi}(G)$; we shall consider it with the operator space structure inherited from $\PM_{p,\pi}(G)^\ast$. It is immediate that $\PM_{p,,\pi}(G)_\ast$ consists of continuous functions on $G$. We have:
\begin{lemma} \label{mullem}
Let $G$ be a locally compact group, let $p \in (1,\infty)$, let $(\pi,E), (\rho, F) \in \Rep_{p'}(G)$, and let $(\pi \tensor \rho, E \ttensor_p F)$ be as in Proposition \emph{\ref{tprop}}. Then pointwise multiplication induces a completely bounded map from $\PM_{p',\pi}(G)_\ast \Tensor \PM_{p',\rho}(G)_\ast$ into $\PM_{p', \pi \tensor \rho}(G)_\ast$ with $\cb$-norm at most $\KG^2$.
\end{lemma}
\begin{proof}
It follows from \cite[Corollary 3.2]{RunBp} that pointwise multiplication of two functions in $\PM_{p',\pi}(G)_\ast$ and $\PM_{p',\rho}(G)_\ast$, respectively, does indeed yield a function in $\PM_{p', \pi \tensor \rho}(G)_\ast$.
\par 
A diagram chase just as in the proof of \cite[Lemma 6.2]{LNR}---invoking Proposition \ref{tprop}(iv)---then shows that the induced bilinear map is indeed completely bounded with norm at most $\KG^2$.
\end{proof}
\par 
We can now prove:
\begin{proposition} \label{mulprop}
Let $G$ be a locally compact group, let $p \in (1,\infty)$, let $(\pi,E), (\rho, F) \in \Rep_{p'}(G)$, and let $(\pi \tensor \rho, E \ttensor_p F)$ be as specified in Proposition \emph{\ref{tprop}}. Then pointwise multiplication induces a completely bounded, bilinear map from $\PF_{p',\pi}(G)^\ast \times \PF_{p',\rho}(G)^\ast$ into $\PF_{p', \pi \tensor \rho}(G)^\ast$ with $\cb$-norm at most $\KG^2$. Moreover, this multiplication is separately continuous with respect to the weak$^\ast$ topologies involved.
\end{proposition}
\begin{proof}
Let
\[
  m \!: \PM_{p',\pi}(G)_\ast \times \PM_{p',\rho}(G)_\ast \to \PM_{p', \pi \tensor \rho}(G)_\ast
\]
denote pointwise multiplication and recall from Lemma \ref{mullem} that $\| m \|_\cb \leq \KG^2$. As a bilinear map between Banach spaces, $m$ has two Arens extensions
\begin{multline*}
  m_1^{\ast\ast} \!:  \PM_{p',\pi}(G)^\ast \times \PM_{p',\rho}(G)^\ast \to \PM_{p', \pi \tensor \rho}(G)^\ast
  \\ \text{and}\qquad 
  m_2^{\ast\ast} \!:  \PM_{p',\pi}(G)^\ast \times \PM_{p',\rho}(G)^\ast \to \PM_{p', \pi \tensor \rho}(G)^\ast.
\end{multline*}
(This construction is mostly done only for the product of a Banach algebra---see \cite{Dal}---, but works as well for general bilinear maps: see \cite{Gro}). It is routinely checked that $m_1^{\ast\ast}$ and $m_2^{\ast\ast}$ are both completely bounded with $\| m_j^{\ast\ast} \|_\cb \leq \KG^2$ for $j=1,2$.
\par 
For $\sigma \in \{ \pi, \rho, \pi \tensor \rho \}$, let $Q_\sigma \!: \PM_{p',\sigma}(G)^\ast \to \PF_{p',\sigma}(G)^\ast$ denote the restriction map, and note that it is a complete quotient map. We claim that
\[
  Q_{\pi \tensor \rho} \circ m_1^{\ast\ast} \!: \PM_{p',\pi}(G)^\ast \times \PM_{p',\rho}(G)^\ast \to \PF_{p', \pi \tensor \rho}(G)^\ast
\]
drops to a map
\[
  \tilde{m} \!: \PF_{p',\pi}(G)^\ast \times \PF_{p',\rho}(G)^\ast \to \PF_{p', \pi \tensor \rho}(G)^\ast,
\]
which is easily seen to be pointwise multiplication and clearly satisfies $\| \tilde{m} \| \leq \KG^\ast$. (We could equally well work with $m_2$.)
\par 
For $\sigma \in \{ \pi, \rho, \pi \tensor \rho \}$, let $\iota_\sigma \!: \PM_{p',\sigma}(G)_\ast \to L_\infty(G)$ and 
$\tilde{Q}_\sigma\!: \PM_{p',\sigma}(G)^\ast \to L_\infty(G) = L_1(G)^\ast$ denote the canonical inclusion and restriction maps, respectively. Also, let $Q \!: L_\infty(G)^{\ast\ast} \to L_\infty(G)$ be the canonical restriction map, and note that it is an algebra homomorphism.
\par 
As the diagram
\[
  \begin{diagram}
  \PM_{p',\pi}(G)_\ast & \times & \PM_{p',\rho}(G)_\ast & \rTo^m & \PM_{p', \pi \tensor \rho}(G)_\ast \\
  \dTo^{\iota_\pi} &        & \dTo_{\iota_\rho} &            & \dTo_{\iota_{\pi \tensor \rho}} \\
  L_\infty(G)          & \times & L_\infty(G)           & \rTo       & L_\infty(G),
  \end{diagram}
\]
where the bottom row is pointwise multiplication in $L_\infty(G)$ commutes, so does
\begin{equation} \label{dia1}
  \begin{diagram}
  \PM_{p',\pi}(G)^\ast & \times & \PM_{p',\rho}(G)^\ast & \rTo^{m_1^{\ast\ast}} & \PM_{p', \pi \tensor \rho}(G)^\ast \\
  \dTo^{Q \circ \iota_\pi^{\ast\ast}} & & \dTo_{Q \circ \iota_\rho^{\ast\ast}} & & \dTo_{Q \circ \iota_{\pi \tensor \rho}^{\ast\ast}} \\
  L_\infty(G)          & \times & L_\infty(G)           & \rTo       & L_\infty(G),
  \end{diagram}
\end{equation}
As 
\[
  Q \circ \iota_\sigma^{\ast\ast} = \tilde{Q}_\sigma \qquad (\sigma \in \{ \pi, \rho, \pi \tensor \rho \}),
\]
this entails the commutativity of 
\[
  \begin{diagram}
  \PM_{p',\pi}(G)^\ast & \times & \PM_{p',\rho}(G)^\ast & \rTo^{m_1^{\ast\ast}} & \PM_{p', \pi \tensor \rho}(G)^\ast \\
  \dTo^{\tilde{Q}_\pi} &        & \dTo_{\tilde{Q}_\rho} &            & \dTo_{\tilde{Q}_{\pi \tensor \rho}} \\
  L_\infty(G)          & \times & L_\infty(G)           & \rTo       & L_\infty(G),
  \end{diagram}
\]
and thus of
\[
  \begin{diagram}
  \PM_{p',\pi}(G)^\ast & \times & \PM_{p',\rho}(G)^\ast & \rTo^{m_1^{\ast\ast}} & \PM_{p', \pi \tensor \rho}(G)^\ast \\
  \dTo^{Q_\pi}         &        & \dTo_{Q_\rho}         &            & \dTo_{Q_{\pi \tensor \rho}} \\
  \PF_{p',\pi}(G)^\ast & \times & \PF_{p',\rho}(G)^\ast & \rTo       & \PF_{p', \pi \tensor \rho}(G)^\ast
  \end{diagram}
\]
with the bottom row being the desired map $\tilde{m}$.
\par 
Finally, since the weak$^\ast$ topology of $\PF_{p',\sigma}(G)^\ast$ for $\sigma \in \{ \pi, \rho, \pi \tensor \rho \}$ coincides with the weak$^\ast$ topology of $L_\infty(G)$ on norm bounded subsets and since multiplication in $L_\infty(G)$ is separately weak$^\ast$ continuous, the commutativity of (\ref{dia1}) and the Kre\u{\i}n--\v{S}mulian theorem (\cite[Theorem V.5.7]{DS}) establish the separate weak$^\ast$ continuity of pointwise multiplication from $\PF_{p',\pi}(G)^\ast \times \PF_{p',\rho}(G)^\ast$ to $\PF_{p', \pi \tensor \rho}(G)^\ast$.
\end{proof}
\par 
Following \cite{RS}, we call a completely bounded Banach algebra a \emph{dual, completely bounded Banach algebra} if it is a dual operator space such that multiplication is separately weak$^\ast$ continuous.
\par 
We can finally state and prove the first theorem of this section:
\begin{theorem} \label{PFthm}
Let $G$ be a locally compact group $G$, and let $p \in (1,\infty)$. Then $\PF_{p'}(G)^\ast$ is a dual, completely bounded Banach algebra with multiplication having the $\cb$-norm at most $\KG^2$.
\end{theorem}
\begin{proof}
By Proposition \ref{mulprop}, pointwise multiplication 
\[
  \PF_{p'}(G)^\ast \times \PF_{p'}(G)^\ast \to \PF_{p',\lambda_{p'} \tensor \lambda_{p'}}(G)^\ast
\]
is completely bounded with $\cb$-norm at most $\KG^2$ and separately weak$^\ast$ continuous.
\par 
From \cite[Theorem 4.6]{LNR} and \cite[Proposition 5.1]{RunBp}, we conclude that $\PF_{p'}(G)$ and $\PF_{p',\lambda_{p'} \tensor \lambda_{p'}}(G)$ are canonically completely isometrically isomorphic. This completes the proof.
\end{proof}
\par 
We shall now turn to the task of turning $B_p(G)$---the $p$-analog of the Fourier--Stieltjes algebra introduced in \cite{RunBp}---into a completely bounded Banach algebra. As in \cite{RunBp}, a difficulty arises due to the fact that $\Rep_{p'}(G)$ is not a set, but only a class; we circumvent the problem, by imposing size restriction on the spaces involved:
\begin{definition}
Let $G$ be a locally compact group, and let $p \in (1,\infty)$. We call $(\pi,E) \in \Rep_{p'}(G)$ \emph{small} if $\mathrm{card}(E) \leq \mathrm{card}(L_1(G))^{\aleph_0}$.
\end{definition}
\begin{remarks}
\item The left regular representation $(\lambda_p,L_p(G))$ is small, as are the cyclic representations used in \cite{RunBp}.
\item Unlike $\Rep_p(G)$, the class of all small representation in $\Rep_p(G)$ is indeed a set.
\end{remarks}
\par 
Let $G$ be a locally compact group, let $p \in (1,\infty)$, and let $(\pi,E), (\rho,F) \in \Rep_p(G)$ be such that $(\rho,F) \subset (\pi,E)$. Then we have a canonical complete contraction from $\PF_{p,\pi}(G)$ to $\PF_{p,\rho}(G)$. Consequently, if $( (\rho_\alpha,F_\alpha))_\alpha$ is a family of representations contained in $(\pi,E)$, we have a canonical complete contraction 
from $\PF_{p,\pi}(G)$ to $\text{$\ell_\infty$-}\bigoplus_\alpha \PF_{p,\rho_\alpha}(G)$.
\par 
We note:
\begin{proposition} \label{small}
Let $G$ be a locally compact group, let $p \in (1,\infty)$, let $(\pi,E) \in \Rep_p(G)$, and let $( (\rho_\alpha,F_\alpha))_\alpha$ be the family of all small representations contained in $(\pi,E)$. Then the canonical map from $\PF_{p,\pi}(G)$ to $\text{$\ell_\infty$-}\bigoplus_\alpha \PF_{p,\rho_\alpha}(G)$ is a complete isometry.
\end{proposition}
\begin{proof}
We need to show the following: for each $n,m \in \posints$, each $n \times n$ matrix $[ f_{j,k} ] \in M_n(L_1(G))$, and each $\epsilon > 0$, there is a closed subspace $F$ of $G$ invariant under $\pi(G)$ with $\mathrm{card}(F) \leq \mathrm{card}(L_1(G))^{\aleph_0}$ such that 
\[
  \left\| \left[ \pi(f_{j,k})^{(m)} |_{M_m(F)} \right] \right\|_{{\cal B}(M_m(F),M_{nm}(F))} \geq 
  \left\| \left[ \pi(f_{j,k})^{(m)} \right] \right\|_{{\cal B}(M_m(E),M_{nm}(E))} - \epsilon.
\]
\par 
Let $n,m \in \posints$, $[ f_{j,k} ] \in M_n(L_1(G))$, and $\epsilon > 0$. Trivially, there is $[ \xi_{\nu,\mu} ] \in M_m(E)$ with
$\| [ \xi_{\nu,\mu} ] \|_{M_m(E)} \leq 1$ such that
\[
  \| [ \pi(f_{j,k})\xi_{\nu,\mu} ] \|_{M_{nm}(E)} \geq 
  \left\| \left[ \pi(f_{j,k})^{(m)} \right] \right\|_{{\cal B}(M_m(E),M_{nm}(E))} - \epsilon.
\]
Let $F$ be the closed linear span of $\{ \pi(f) \xi_{\nu,\mu} : f \in L_1(G), \, \nu,\mu = 1, \ldots, m \}$; it clearly has the desired properties.
\end{proof}
\begin{definition} \label{unidef}
Let $G$ be a locally compact group, and let $p \in (1,\infty)$. We say that $(\pi_u,E_u) \in \Rep_p(G)$ \emph{$p$-universal} if it contains every small reperesentation in $\Rep_p(G)$. We write $\UPF_p(G)$ instead of $\PF_{p,\pi_u}(G)$ and call the elements of $\UPF_p(G)$ \emph{universal $p$-pseudofunctions}.
\end{definition}
\begin{remarks}
\item Since cyclic representations in the sense of \cite{RunBp} are small, a $p$-universal representation according to Definition \ref{unidef} is also $p$-universal in the sense of \cite[Definition 4.5]{RunBp}. We do not know if the converse is also true.
\item There are indeed $p$-universal representations: this can be seen as in the example immediately after \cite[Definition 4.5]{RunBp}.
\item If $(\pi_u,E_u) \in \Rep_p(G)$ is $p$-universal and $(\rho,F) \in \Rep_p(G)$ is arbitrary, then Proposition \ref{small} shows that we have a canonical complete contraction from $\UPF_p(G)$ to $\PF_{p,\rho}(G)$. In particular, the operator space structure of $\UPF_p(G)$ does not depend on a particular $p$-universal representation.
\end{remarks}
\par 
Let $G$ be a locally compact group, and let $p \in (1,\infty)$. As every $p'$-universal representation of $G$ is also $p'$-universal in the sense of \cite{RunBp}, \cite[Theorem 6.6(ii)]{RunBp} remains valid, and we can identity $B_p(G)$ with the Banach space dual of $\UPF_{p'}(G)$. As $\UPF_{p'}(G)$ is an operator space by virtue of Definition \ref{osdef}, we define the canonical operator space structure of $B_p(G)$ as the one it has as the dual space of $\UPF_{p'}(G)$.
\begin{theorem} \label{Bpthm}
Let $G$ be a locally compact group, and let $p \in (1,\infty)$. Then:
\begin{items}
\item $B_p(G)$ is a dual, completely bounded Banach algebra;
\item the canonical image of $\PF_{p'}(G)^\ast$ in $B_p(G)$ is an ideal of $B_p(G)$.
\end{items}
\end{theorem}
\begin{proof}
Let $(\pi_u,E_u) \in \Rep_{p'}(G)$ be $p'$-universal. By Proposition \ref{mulprop}, pointwise multiplication 
\begin{equation} \label{muleq}
  B_p(G) \times B_p(G) \to \PF_{p',\pi_u \tensor \pi_u}(G)^\ast
\end{equation}
is completely bounded. Since $(\pi_u,E_u)$ is $p'$-universal, we have a canonical complete contraction from $\UPF_{p'}(G)$ to $\PF_{p',\pi_u \tensor \pi_u}(G)$. Composing the adjoint of this map with (\ref{muleq}), we obtain pointwise multiplication on $B_p(G)$, which is thus completely bounded. That multiplication in $B_p(G)$ is separately weak$^\ast$ continuous is seen as in the proof of Theorem \ref{PFthm}. This proves (i).
\par 
(ii) follows from \cite[Proposition 5.1]{RunBp}.
\end{proof}
\begin{remarks}
\item Unless $p =2$, we it is well possible that the canonical map from $\PF_{p'}(G)^\ast$ to $B_p(G)$ fails to be an isometry: see the remark immediately after \cite[Corollary 5.3]{RunBp}.
\item We even have that $\PF_{p'}(G)^\ast$ is a $B_p(G)$ module with completely bounded module actions. Since the canonical map from $\PF_{p'}(G)^\ast$ to $B_p(G)$ need not be a (complete) isometry, this is somewhat stronger than Theorem \ref{Bpthm}(ii).
\end{remarks}
\section{Herz--Schur and completely bounded multipliers of $A_p(G)$}
Let $G$ be a locally compact group, let $p,q \in (1,\infty)$, and---as in \cite{ER} and \cite{LNR}---let $\tensor^\gamma$ stand for the projective tensor product of Banach spaces. Even though $L^p(G) \tensor^\gamma L^q(G)$ does not consist of functions on $G \times G$, but rather of equivalence classes of functions, it still makes sense to speak of multipliers of $L^p(G) \tensor^\gamma L^q(G)$: by a multiplier of $L^p(G) \tensor^\gamma L^q(G)$, we mean a continuous function $f$ on $G \times G$, so such that the corresponding multiplication operator $M_f$ induces a bounded linear operator on $L^p(G) \tensor^\gamma L^q(G)$.
\par 
For $p \in (1,\infty)$, we write ${\cal V}_p(G)$ to denote the pointwise multipliers of $L^p(G) \tensor^\gamma L^{p'}(G)$.
For any function $f \!: G \to \comps$, we write $K(f)$ for the function
\[
  G \times G \to \comps, \quad (x,y) \mapsto f(xy^{-1})
\]
We define the \emph{Herz--Schur multipliers} of $A_p(G)$ as
\[
  \MHS(A_p(G)) := \{ f \!: G \to \comps : K(f) \in {\cal V}_p(G) \}.
\]
As ${\cal V}_p(G)$ is a closed subspace of ${\cal B}(L^p(G) \tensor^\gamma L^{p'}(G))$, and since the map $\MHS(A_p(G)) \ni f \mapsto K(f)$ is injective, we can equip $\MHS(A_p(G))$ with a natural norm turning it into a Banach space.
\par
In \cite{BF}, M.\ Bo\.{z}ejko and G.\ Fendler showed that $\MHS(A(G))$ and $\Mcb(A(G))$ are isometrically isomorphic (see also \cite{Jol}), and in \cite{Fen}, Fendler showed that, for general $p \in (1,\infty)$, the Herz--Schur multipliers of $A_p(G)$ are precisely the $p$-completely bounded ones.
\par 
In this section, we investigate how $\MHS(A_p(G))$ and $\Mcb(A_p(G))$ relate to each other for general $p \in (1,\infty)$, but with $A_p(G)$ having the operator space structure introduced in \cite{LNR}. We start with a lemma:
\begin{lemma} \label{HSlem1}
Let $p \in (1,\infty)$, let $X$ and $Y$ be measure spaces, and let $E$ be a $\QSL_p$-space. Then the map
\begin{multline*}
  (L^p(X) \tensor E) \tensor (L^{p'}(Y) \tensor E^\ast) \to L^p(X) \tensor L^{p'}(Y), \\
  (f \tensor \xi) \tensor (g \tensor \phi) \mapsto \langle \xi, \phi \rangle (f \tensor g)
\end{multline*}
extends to a complete quotient map 
\[
  \tr_E \!: \ROW(L^p(X,E)) \Tensor \COL(L^{p'}(Y,E^\ast)) \to \ROW(L^p(X)) \Tensor \COL(L^{p'}(Y)).
\]
\end{lemma}
\begin{proof}
Since $E^\ast$ is a $\QSL_{p'}$-space, this follows from \cite[Theorem 4.6]{LNR} through taking adjoints.
\end{proof}
\begin{proposition} \label{HSprop}
Let $p \in (1,\infty)$, and let $G$ be a locally compact group. Then, for any $f \in \MHS(A_p(G))$, the multiplication operator $M_{K(f)} \!: L^p(G) \tensor^\gamma L^{p'}(G) \to L^p(G) \tensor^\gamma L^{p'}(G)$ is completely bounded on $\ROW(L^p(G)) \Tensor \COL(L^{p'}(G))$ such that
\[
  \| M_{K(f)} \|_\cb   = \| f \|_{\MHS(A_p(G))}.
\]
\end{proposition}
\begin{proof}
Let $f \in \MHS(A_p(G))$, and let $\epsilon > 0$. Then, by \cite{Gil} (see also \cite[Theorem 4.4]{Fen}), there is a $\QSL_p$-space $E$ along with bounded continuous functions $L \!: G \to E$ and $R \!: G \to E^\ast$ such that
\[
  K(f)(x,y) = \langle L(x), R(y)\rangle \qquad (x,y \in G)
\]
and
\begin{equation} \label{estim1}
  \| L \|_\infty \| R \|_\infty < \| f \|_{\MHS(A_p(G))} + \epsilon,
\end{equation}
where 
\[
  \| L \|_\infty := \sup_{x \in G} \| L(x) \| \qquad\text{and}\qquad \| R \|_\infty := \sup_{x \in G} \| R(x) \|.
\]
\par
Define $\tilde{L} \!: L^p(G) \to L^p(G,E)$ by letting
\[
  (\tilde{L}\xi)(x) := \xi(x)L(x) \qquad (\xi \in L^p(G), \, x \in G).
\]
Then $\tilde{L}$ is linear and bounded with $\| \tilde{L} \| = \| L \|_\infty$. Similarly, one defines a bounded linear map $\tilde{R} \!: L^{p'}(G) \to L^{p'}(G,E^\ast)$ with $\| \tilde{R} \| = \| R \|_\infty$ by letting
\[
  (\tilde{R}\eta)(x) := \eta(x)R(x) \qquad (\eta \in L^{p'}(G), \, x \in G).
\]
Since the row and the column spaces over any Banach space are homogeneous, it is clear that $\tilde{L} \!: \ROW(L^p(G)) \to \ROW(L^p(G,E))$ and $\tilde{R} \!: \COL(L^{p'}(G)) \to \COL(L^{p'}(G,E^\ast))$ are completely bounded with
$\| \tilde{L} \|_\cb = \| L \|_\infty$ and $\| \tilde{R} \|_\cb = \| R \|_\infty$. From \cite[Corollary 7.1.3]{ER}, it thus follows that 
\[
  \tilde{L} \tensor \tilde{R} \!:   \ROW(L^p(G)) \Tensor \COL(L^{p'}(G)) \to 
  \ROW(L^p(G,E)) \Tensor \COL(L^{p'}(G,E^\ast)) 
\]
is completely bounded as well with $\| \tilde{L} \tensor \tilde{R} \|_\cb \leq \| L \|_\infty \| R \|_\infty$. Since 
\[
  \tr_E \!: \ROW(L^p(G,E)) \Tensor \COL(L^{p'}(G,E^\ast)) \to   \ROW(L^p(G)) \Tensor \COL(L^{p'}(G))
\]
as in Lemma \ref{HSlem1} is a complete contraction, we conclude that $\tr_E \circ (\tilde{L} \tensor \tilde{R})$ is completely bounded with $\cb$-norm at most $\| L \|_\infty \| R \|_\infty$. 
\par
From the definitions of $\tr_E$, $\tilde{L}$, and $\tilde{R}$, it it straightforward to verify that $\tr_E \circ (\tilde{L} \tensor \tilde{R}) = M_{K(f)}$. In view of (\ref{estim1}), we thus obtain that
\[
  \| M_{K(f)} \|_\cb \leq \| L \|_\infty \| R \|_\infty <  \| f \|_{\MHS(A_p(G))} + \epsilon.
\]
Since $\epsilon > 0$ is arbitrary, this yields $\| M_{K(f)} \|_\cb \leq \|  f  \|_{\MHS(A_p(G))}$. By definition, 
\[
  \| f \|_{\MHS(A_p(G))} = \| M_{K(f)} \| \leq \| M_{K(f)} \|_\cb,
\]
holds, so that have equality as claimed.
\end{proof}
\par
Passing to quotients we thus obtain:
\begin{theorem} \label{HSthm}
Let $p \in (1,\infty)$, and let $G$ be a locally compact group. Then $\MHS(A_p(G))$ is contained in $\Mcb(A_p(G))$ such that the inclusion is a contraction.
\end{theorem}
\begin{remarks}
\item By Proposition \ref{HSprop}, the linear map $\MHS(A_p(G)) \ni f \mapsto M_{K(f)}$ is an isometric embedding into the operator space $\CB(\ROW(L^p(G)) \Tensor \COL(L^{p'}(G)))$ and can be used to equip $\MHS(A_p(G))$ with a canonical operator space structure. We do not know whether Theorem \ref{HSthm} can be improved to yield a completely contractive---or at least completely bounded---inclusion map.
\item For amenable $G$, the algebras $\PF_{p'}(G)^\ast$, $B_p(G)$, $\Mcb(A_p(G))$, and $\MHS(A_p(G))$ are easily seen to be isometrically isomorphic (see \cite{Cow}, \cite{Her3}, and \cite{RunBp}). We do not know whether theses isometric isomorphisms are, in fact, completely isometric; for some of them, this seems to be open even in the case where $p =2$.
\end{remarks}
\par 
For $p \in (1,\infty)$ and a locally compact group $G$, any $f \in A_p(G)$ is a $\cb$-multiplier of $A_p(G)$, simply because $A_p(G)$ is a completely bounded Banach algebra. However, as $A_p(G)$ is not known to be completely contractive, this does not allow us
to conclude that $\|f \|_{\Mcb(A_p(G))} \leq \| f \|_{A_p(G)}$, but only that $\|f \|_{\Mcb(A_p(G))} \leq \KG^2 \| f \|_{A_p(G)}$.
\par
Theorem \ref{HSthm}, nevertheless allows us to obtain a better norm estimate:
\begin{corollary} 
Let $p \in (1,\infty)$, and let $G$ be a locally compact group. Then we have
\[
  \| f \|_{\Mcb(A_p(G))} \leq \| f \|_{A_p(G)} \qquad (f \in A_p(G)).
\]
\end{corollary}
\begin{proof}
Let $f \in A_p(G)$. By \cite[Proposition 10.2]{Pie}, we have $\| f \|_{\MHS(A_p(G))} \leq \| f \|_{A_p(G)}$ and thus
\[
  \| f \|_{\Mcb(A_p(G))} \leq \| f \|_{\MHS(A_p(G))} \leq \| f \|_{A_p(G)}
\]
by Theorem \ref{HSthm}.
\end{proof}
\section{$B_p(G)$, $\PF_{p'}(G)^\ast$, and the amenability of $G$}
A classical amenability criterion due to R.\ Godement asserts that a locally compact group $G$ is amenable if and only if its trivial representation is weakly contained in $(\lambda_2,L_2(G))$ (see \cite[Theorem 8.9]{Pie}, for instance). In terms of Fourier--Stieltjes algebras this means that $G$ is amenable if and only if $B_r(G) = B(G)$ (the equality is automatically an complete isometry).
\par 
The following theorem generalizes this to a general $L_p$-context:
\begin{theorem} \label{amthm}
The following are equivalent for a locally compact group $G$:
\begin{items}
\item the canonical map from $\PF_{p'}(G)^\ast$ into $B_p(G)$ is surjective for each $p \in (1,\infty)$;
\item there is $p \in (1,\infty)$ such that $1 \in \PF_{p'}(G)^\ast$;
\item $G$ is amenable.
\end{items}
\end{theorem}
\begin{proof}
(i) $\Longrightarrow$ (ii) is trivial.
\par 
(ii) $\Longrightarrow$ (iii): An inspection of the proof of \cite[Theorem 5]{Cow} shows that $1 \in \PF_{p'}(G)^\ast$ for just one $p \in (1,\infty)$ is possible only if $G$ is amenable.
\par 
(iii) $\Longrightarrow$ (i): This follows from \cite[Theorem 6.7]{RunBp}.
\end{proof}
\begin{remark}
Unless $p=2$, we cannot say for amenable $G$ whether or not $\PM_{p'}(G)^\ast = B_p(G)$ holds completely isometrically. By \cite[Theorem 6.7]{RunBp}, we do have an isometric isomorphism, but this is all we can say. Due to the lack of an inverse mapping theorem for completely bounded maps, we even do not know for general $p \in (1,\infty)$ whether the completely bounded bijective map from $\PF_{p'}(G)^\ast$ onto $B_p(G)$ has a completely bounded inverse.
\end{remark}
\par 
In \cite{Rua}, Z.-J.\ Ruan adapted the notion of an amenable Banach algebra due to B.\ E.\ Johnson (\cite{Joh}) to an operator space context. 
\par 
Given a completely bounded Banach algebra $\A$ and a completely bounded Banach $\A$-bimodule $E$, i.e., a Banach $\A$-bimodule which is also an operator space such that the module actions are completely bounded, the dual operator space $E^\ast$ because a completely bounded Banach $\A$-bimodule in its own right via
\[
  \langle \xi, a \cdot \phi \rangle := \langle \xi \cdot a, \phi \rangle 
  \quad\text{and}\quad
  \langle \xi, \phi \cdot a \rangle := \langle a \cdot \xi, \phi \rangle 
  \qquad (\xi \in E, \, \phi \in E^\ast, \, a \in \A),
\] 
and $\A$ is said to be \emph{operator amenable} if and only if, for each completely bounded Banach $\A$-bimodule $E$, every completely bounded derivation $D \!: \A \to E^\ast$ is inner. Ruan showed that a locally compact group $G$ is amenable if and only if $A(G)$ is operator amenable, and in \cite{LNR}, Lambert and the authors extended this result to $A_p(G)$ for arbitrary $p \in (1,\infty)$.
\par 
Suppose that $\A$ is a \emph{dual}, completely bounded Banach algebra. If $E$ is a completely bounded Banach $\A$-bimodule, we call $E^\ast$ \emph{normal} if the bilinear maps 
\[
  \A \times E^\ast \to E^\ast, \quad \quad (a,\phi) \mapsto \left\{ \begin{array}{c} a \cdot \phi, \\ \phi \cdot a \end{array}
  \right.
\] 
are separately weak$^\ast$ continuous. Following \cite{RS}, we say that $\A$ is \emph{operator Connes-amenable} if, for every completely bounded Banach $\A$-bimodule $E$ such that $E^\ast$ is normal, every weak$^\ast$-weak$^\ast$-continuous, completely bounded derivation $D \!: \A \to E^\ast$ is inner.
\par 
Extending \cite[Theorem 4.4]{RS} in analogy with \cite[Theorem 7.3]{LNR}, we obtain:
\begin{theorem}
The following are equivalent for a locally compact group $G$:
\begin{items}
\item $G$ is amenable;
\item $P_{p'}(G)^\ast$ is operator Connes-amenable for every $p \in (1,\infty)$;
\item $B_r(G)$ is operator Connes-amenable;
\item there is $p \in (1,\infty)$ such that $\PF_{p'}(G)^\ast$ is operator Connes-amenable.
\end{items}
\end{theorem}
\begin{proof}
(i) $\Longrightarrow$ (ii): Let $p \in (1,\infty)$ be arbitrary. Then \cite[Theorem 7.3]{LNR} yields the operator amenability of $A_p(G)$. Since the inclusion of $A_p(G)$ into $\PF_{p'}(G)^\ast$ is a completely contractive algebra homomorphism with weak$^\ast$ dense range, the operator space analog of \cite[Proposition 4.2(i)]{RunD} yields the operator Connes-amenability of $\PF_{p'}(G)^\ast$.
\par 
(ii) $\Longrightarrow$ (iii) $\Longrightarrow$ (iv) are trivial.
\par 
(iv) $\Longrightarrow$ (i): Let $p \in (1,\infty)$ be such that $\PF_{p'}(G)^\ast$ is operator Connes-amenable. The operator space analog of \cite[Proposition 4.1]{RunD} then yields that $\PF_{p'}(G)^\ast$ has an identity, so that Theorem \ref{amthm}(ii) is satisfied. By Theorem \ref{amthm}, this means that $G$ is amenable.
\end{proof}
\renewcommand{\baselinestretch}{1.0}
\dated
\renewcommand{\baselinestretch}{1.2}
\vfill
\begin{tabbing} 
\textit{Second author's address}: \= Department of Mathematical and Statistical Sciences \kill
\textit{First author's address}:  \> School of Mathematics and Statistics \\
                                  \> 4364 Herzberg Laboratories \\ 
                                  \> Carleton University \\ 
                                  \> Ottawa, Ontario \\
                                  \> Canada K1S 5B6 \\[\medskipamount]
\textit{E-mail}:                  \> \texttt{mneufang@math.carleton.ca} \\[\bigskipamount]  
\textit{Second author's address}: \> Department of Mathematical and Statistical Sciences \\
                                  \> University of Alberta \\
                                  \> Edmonton, Alberta \\
                                  \> Canada T6G 2G1 \\[\medskipamount]
\textit{E-mail}:                  \> \texttt{vrunde@ualberta.ca}
\end{tabbing}
\end{document}